\documentclass{amsart}[11pt]
\usepackage{mathrsfs,amssymb,amsfonts,amsmath}
\newtheorem{teo}{Theorem}[section]

\newtheorem{lem}{Lemma}[section]
\newtheorem{nota}{Remark}[section]
\newcommand{\re}{\mathbb{R}}
\title{Generating functions of Cauchy-Stieltjes type for orthogonal polynomials}
\date{\today}

\begin{document}
\maketitle
\centerline{\author{Marek Bozejko\footnote{Instytut Matematyczny, Uniwersytet Wroclawski, Pl. Grunwaldzki 2/4, 50-384 Wroclaw, Poland.} and
 Nizar Demni \footnote{Faculty of Mathematics, Bielefeld University, Germany.\\ 
 The first author was partially supported by KBN grant no 1 PO3A 01330 and by a Marie Curie Transfer of Knowledge Fellowship of the European Community 's Sixth Framework porgram under contract number MTKD-CT-2004-013389}}}

\begin{abstract}
We give a free probabilistic interpretation of the multiplicative renormalization method applied with $h(x) = (1-x)^{-1}$. As a byproduct, we give a short proof of the  Asai-Kubo-Kuo result on the characterization of the family of probability measures with a Cauchy-Stieltjes type generating functions for orthogonal polynomials. This family turns to be the free Meixner family. We also give representations for the Voiculescu transform of all free Meixner distributions (even in the non-freely infinitely divisible case).    
\end{abstract}

\section{Multiplicative renormalization method}
The theory of orthogonal polynomials has seen a significant  growth during the last century due to their connections with other areas such as approximation theory, quantum probability, stochastic processes (see \cite{Acc}, \cite{Sch}). Nevertheless, for a given measure, the computation of the corresponding family of orthogonal polynomials remains a quite hard and often even impossible task by the use of Gram-Schmidt orthogonalization method (consider for instance the geometric distribution). This was the starting point and the motivation of some reserachers to find out a more subtle method to derive them. In \cite{Asai2}, the authors propose  the so-called \emph{multiplicative renormalization} method based on generating functions. More precisely, let $\mu$ an infinitely supported measure on the real line with finite all order moments. Then, the family $(P_n)_{n \geq 0}$ of orthogonal \emph{monic} polynomials with respect to $\mu$ satisfies the three terms-recurrence relation 
\begin{equation*}
xP_n(x) = P_{n+1}(x) + \alpha_nP_n(x) + \omega_nP_{n-1}(x)
\end{equation*} 
where $\alpha_n \in \re, \omega_0 = 1, \omega_n > 0$ for all  $n\geq 1$. $(\alpha_n)_n, (\omega_n)_n$ are called the Jacobi-Szeg\"o parameters. $\mu$ is symmetric if and only if $\alpha_n = 0$ for all $n$. A {\it pre-generating function} for $\mu$ is an infinite series of the form
\begin{equation*}
\phi(z,x) := \sum_{n \geq 0}S_n(x)z^n
\end{equation*}
where 
\begin{itemize}
\item $S_n(x)$ is a polynomial of degree $n$ for each $n \geq 0$.
\item $\limsup_{n \rightarrow \infty} ||S_n||_{L^2(\mu)} < \infty$.
\end{itemize}
When $(S_n)_{n \geq 0}$ are orthogonal with respect to $\mu$, then $(z,x) \mapsto \phi(z,x)$ is said to be a {\it generating function}. Otherwise, its {\it multilplicative renormalization}  is defined by
\begin{equation*}
\psi(z,x) := \frac{\phi(z,x)}{\mathbb{E}(\phi(z,X))} = \sum_{n \geq 0}Q_n(x)z^n
\end{equation*}
where $X$ is a random variable in some probability space with law $\mu$ and $\mathbb{E}$ denotes the expectation. 
In order to decide  whether a given function $\psi$ is a generating function for a pre-given measure $\mu$, Asai-Kuo-Kubo (\cite{Asai1},\cite{Asai1b},\cite{Asai2},\cite{Asai3}) derived the following criterion:
\begin{teo}\label{AKK}
$\psi$ is a generating function for $\mu$ if and only if $\mathbb{E}(\psi(z,X)\psi(v,X))$ only depends on $zv$. 
\end{teo}
In that case, one writes $Q_n(x) = a_nP_n(x)$, the coefficients $(a_n)_n$ being determined by 
\begin{equation}
\label{coeff}
\lim_{z \rightarrow 0} \psi(z,\frac{x}{z}) = \sum_{n \geq 0}a_nx^n
\end{equation}

When considering most of classical orthogonal polynomials, the authors distinguished:
\begin{equation*}
\phi(z,x) = e^{\rho(z)x},\quad \phi(z,x) = (1- \rho(z)x)^{-\beta},\, \beta > 0
\end{equation*}
where $\rho$ is an analytic function around $0$ with $\rho(0) = 0$ and $\rho'(0) \neq 0$. $\phi$ is of the form $h(\rho(z)x)$ for some function $h$ so that when Theorem \ref{AKK} holds, we say that the multiplicative renormalization method applies with $h$.  In \cite{Kubo1},\cite{Kubo2}, the authors answered the following question: characterize the family of probability measures applicable with $h(x) = (1-x)^{-1}$. The Jacobi-Szeg\"o parameters $(\alpha_n)_n, (\omega_n)_n$ are shown to be stationary sequences from rank $n=0, n=1$ respectively, a fact that characterizes the so-called free Meixner distributions. These distributions form the so-called free Meixner family defined in \cite{Ans} via their generating functions and they occured under other different guises. They are for instance solutions of a quadratic regression problem in free probability 
\cite{Boz} (see also \cite{Bry} for another characterization). As a matter of fact, we shall use free probability theory to explain the occurrence of the free Meixner family in \cite{Kubo1}, 
\cite{Kubo2}. The techniques we use not only match \cite{Kubo1}, \cite{Kubo2} to \cite{Ans}, but also give an elegant and easy proof for the representations of the Voiculescu transforms for the free Meixner distributions. These representations was already derived in \cite{Bo} only for freely infinitely divisible laws from the free Meixner family while here, we deal with all the distributions.\\ 
The outline of the paper is as follows. We show how the framework in \cite{Kubo1} and in \cite{Kubo2} fits into the one in \cite{Ans} using free probability theory. This leads as a byproduct to the free Meixner family. Finally, we derive using Asai-Kuo-Kubo criterion the representationa of their Voiculsecu transforms. 
\begin{nota}
When $h(x) = e^x$, the use of Theorem \ref{AKK} gives a representation of  the classical cumulant generating function which was already proved for classical Meixner laws (\cite{Kubo}, \cite{Sch}). The case $h(x) = (1-x)^{-\beta}, \beta > 0$ is dealt with in a forthcoming paper by the second author (\cite{Dem1}). We also hint the interested reader to the survey-like paper by the authors (\cite{Bozejko}) where helpful remarks on the the classical, free and quantum Meixner families are given.    
\end{nota}

\section{Free Meixner family and multiplicative renormalization method}
In this section, we formulate the multiplicative renormalization method applicable with $h(x) = (1-x)^{-1}$ using Cauchy-Stieltjes transforms. Once we did, the free Meixner family, recalled below, will follow by the simple application of Lemma 2 p. 243 in \cite{Ans}. To proceed, some needed facts and definitions from probability theory are collected in the following paragraph.
\subsection{Definitions and needed facts}  The Cauchy-Stieltjes or $G$-transform of a measure $\mu$ is defined by
\begin{equation*}
G_{\mu}(z) := \int_{\re}\frac{1}{z-x} \mu(dx), z \in \mathbb{C} \setminus \textrm{supp}(\mu).
\end{equation*}
It maps  the open upper half plane $\mathbb{C}_+$ into the open lower half plane $\mathbb{C}_-$ and it is analytic there. 
Let $F_{\mu} := 1/G_{\mu}$, then it was shown in \cite{Ber1}, \cite{Ber2} that $F_{\mu}$ is one-to-one in some neighborhood of infinity. The Voiculescu transform is then defined as 
\begin{equation*}
\phi_{\mu}(z) = F_{\mu}^{-1}(z) - z 
\end{equation*} 
where $F_{\mu}^{-1}$ is the right inverse of $F_{\mu}$. Equivalently, one may use the \emph{free-cumulant generating function} defined by
\begin{equation*}
R_{\mu}(z) : = \phi_{\mu}(1/z) = K_{\mu}(z) - \frac{1}{z}, \quad K_{\mu}(z) := F_{\mu}^{-1}(1/z),
\end{equation*}

Suppose that $\mu$ is probability measure with all-order finite moments, its continued fraction expansion is written (see \cite{Hora})
\begin{equation}\label{CF}
G_{\mu}(z) = \frac{1}{z - \alpha_0 - \displaystyle - \frac{w_1}{z - \alpha_1 - \displaystyle \frac{\omega_2}{z- \alpha_2 - \displaystyle \frac{\omega_3}{z-\dots}}}}.
\end{equation}
\subsection{Free Meixner family}
The free Meixner family  was defined in \cite{Ans} and can be characterized using Jacobi-Szeg\"o parameters. Otherwise stated, we shall assume from now on that $\mu$ is a standard probability measure, that is with zero mean ($\alpha_0 = 0$) and unit variance ($\omega_1=1$). Adotpting the notations used in \cite{Bo}, $\mu$ belongs to the free Meixner family if and only if its Jacobi-Szeg\"o parameters are given by
\begin{equation*}
\alpha_n = a \in \re,a \in \re, n \geq 1,\, \omega_n = 1+b, b \geq -1, n \geq 2.
\end{equation*} 
$\mu$ is then referred to as a free Meixner distribution of parameters $a,b$. The value $b=-1$ corresponds to a Dirac mass. There are, up to affine transformations, six compactly-supported  probability distributions, only five of them are freely infinitely divisible (\cite{Sait}).   

\subsection{Multiplicative renormalization method and free probability}
Choosing $h = (1-x)^{-1}$, Theorem \ref{AKK} reads 
\begin{teo}\label{AKK1}
Let $\mu$ be a standard probability measure with all-order finite moments. Let $\rho$ be an analytic function around zero with $\rho(0) = 0, \rho'(0) \neq 0$ and $X$ be a random variable with law $\mu$. Then 
\begin{equation*}
\psi(z,x) = \frac{(1-\rho(z)x)^{-1}}{\mathbb{E}[1-\rho(z)X)^{-1}]}
\end{equation*}
is a generating function for $\mu$ if and only if for all $z,v$ in a neighborhood of zero 
\begin{equation*}
\frac{\mathbb{E}[[(1-\rho(z)X)(1-\rho(v)X)]^{-1}]}{\mathbb{E}[1-\rho(z)X)^{-1}]\mathbb{E}[1-\rho(v)X)^{-1}]}
\end{equation*}
depends only on $zv$.
\end{teo}

In order to translate the above Theorem to free probability theory, we first define 
\begin{equation*}
1/g(z) := F_{\mu}(1/\rho(z)) \quad \Leftrightarrow \quad 1/\rho(z) = F_{\mu}^{-1}(1/g(z)) = K_{\mu}(g(z)).
\end{equation*}
Then, we recall the relation (\cite{Ber1},\cite{Ber2})
\begin{equation*}
F_{\mu}(z) = z- \phi_{\mu}(F_{\mu}(z)) \quad (\textrm{or} \quad G_{\mu}(z) = \frac{1}{z- R_{\mu}(G_{\mu}(z))}),
\end{equation*}
and the elemantary transformation
\begin{equation*}
\int \frac{1}{1-zx}\mu(dx) = \frac{1}{z}G_{\mu}\left(\frac{1}{z}\right). 
\end{equation*}
As a result, after elemantary computations, Theorem \ref{AKK1} translates to
\begin{teo}
\label{BD}
Let $\mu$ be a standard probability measure with all-order finite moments. Let $g$ be analytic around zero with $g(0) = 0$ and $g'(0) \neq 0$. Then  
\begin{equation*}
\frac{1}{g(z)(K_{\mu}(g(z)) - x)}
\end{equation*}
is a generating function for $\mu$ if and only if
\begin{equation*}
\frac{\phi_{\mu}(1/g(z)) - \phi_{\mu}(1/g(v))}{1/g(z) - 1/g(v)}
\end{equation*}
only depends on $zv$. 
\end{teo}

Finally, one may assume that $\rho'(0) = 1$ or equivalently $g'(0) = 1$ so that 
\begin{equation*}
\psi(z,x) = \sum_{n\geq 0}P_n(x)z^n
\end{equation*}
and use Lemma 2 p. 243 in \cite{Ans} that we recall below: 
\begin{lem}
Let $\mu$ be a (non-necessarily standard) probability measure whose all moments are finite. Let $(P_n)_n$ be the monic polynomials orthogonal with respect to $\mu$. Then,
\begin{equation*}
\sum_{n \geq 0} P_n(x) z^n = \frac{1}{u(z)(f(z) - x)}
\end{equation*}
for $u,uf$ having formal power series expansion with $u(0) = 0, u'(0) = 1, [uf](0) = 1$ if and only if 
\begin{equation*}
\alpha_n = \alpha - \delta_{n0}\alpha', n \geq 0, \omega_n = \beta - \delta_{n1}\beta', n \geq 1,
 \end{equation*}
for some $\alpha,\alpha',\beta,\beta'$ with $\beta-\beta' \geq 0$. Moreover, $f(z) = K_{\mu}(u(z))$.
\end{lem}
Since the Jacobi-Szeg\"o parameters characterize $\mu$ and identifying $a,b$ with $\alpha',\beta'$ when $\mu$ is standard, we finally get: 
\begin{teo}
Let $\mu$ be a probability measure whose all moments are finite. Then, the multiplicative renormalization method applies for $\mu$ with $h(x) = (1-x)^{-1}$ if and only if $\mu$ is a free Meixner distribution.  
\end{teo}

\section{Representations of Voiculescu transforms}
Using results in \cite{Sait}, authors in \cite{Bo} showed that free Meixner distributions of parameters $a\in \re,b \geq 0$ are freely infinitely divisible. Moreover, the L\'evy-Kinchin measure is identified to a semi-circle law of mean $a$ and variance $b$. This  yielded representations of their Voiculescu transforms as Cauchy-Stieltjes transforms. 
Here, we use Theorem \ref{AKK1} to recover these representations and to investigate the remaining, non freely infinitely, case. \\
From the proof of Lemma 2. in \cite{Ans}, one has (with the same notations used there)
\begin{equation*}
u(z) = \frac{z}{1+ \alpha'z + \beta'z^2} 
\end{equation*}
or using the parameters $a,b$ that 
\begin{equation*}
g(z) = \frac{z}{1+ az + bz^2} .
\end{equation*}
It follows from the second assertion in Theorem \ref{BD} that
\begin{equation*}
\frac{\phi_{\mu}((bz^2 + az + 1)/z) - \phi_{\mu}((bv^2 + av + 1)/v)}{z-v}
\end{equation*}
depends only on $zv$. Taking $v = 0$, since $\phi_{\mu}(\infty) = R_{\mu}(0)  = 0$ (mean), then 
\begin{equation*}
\phi_{\mu}((bz^2 + az + 1)/z) = Cz, \quad z \in V(0).
\end{equation*}
Writing $\phi_{\mu}(z) = R_{\mu}(1/z)$ and since $R_{\mu}(y)/y \rightarrow  1$ (variance) as $y \rightarrow 0$, then it follows by letting $z \rightarrow 0$ that $C=1$.  
The obtained relation is two-fold. \\On the one hand, $\mu$ is identified through its continued fraction as follows: let $w = (b+1)z + a + 1/z$, then 
\begin{equation*}
G_{\mu}(w) = \frac{1}{F_{\mu}(w)} = \frac{1}{bz + a + 1/z} = \frac{1}{w -z}.
\end{equation*}
Using successively $z = 1/(w - a - (1+b)z)$, one gets the continued fraction 
\begin{equation*}
G_{\mu}(w) = \frac{1}{w - \displaystyle \frac{1}{w- a - \displaystyle \frac{1+b}{w - \dots}}}
\end{equation*}
so that $\mu$ is a standard free Meixner distribution. On the other hand, let $y= 1/g(z)$, then for $b \geq 0$,
\begin{equation*}
z  = \frac{1}{y-a -bz} = \frac{1}{y-a - \displaystyle \frac{b}{y - a - bz}} =\dots =  G_{\omega_{a,b}}(y).
\end{equation*}
Since 
\begin{equation*}
z = \phi_{\mu}(1/g(z)) = \phi_{\mu}(y)
\end{equation*} 
we get the desired representation of the Voiculescu transforms:
\begin{equation*}
 \phi_{\mu}(y) = G_{\omega_{a,b}}(y).
 \end{equation*}
When $-1 \leq b < 0$, then  
\begin{equation*}
K_{\mu}(g(z)) = \frac{(1+b)z^2 + az + 1}{z} = K_{\omega_{a,1+b}}(z).
\end{equation*}
It follows that 
\begin{equation*}
g^{-1}[G_{\mu}(z)] = G_{\omega_{a,1+b}}(z) = \phi_{\eta}(z) \Leftrightarrow G_{\mu}(z) = g[\phi_{\eta}(z)],
\end{equation*} 
in some neighborhood of $\infty$, where $\eta$ is a free Meixner distribution of parameters $a \in \re,1+b \in ]0,1[$.

{\bf Acknowledgments.} Authors thanks Professor F. G\"otze for the fruitful  work atmosphere and the CRC 701  for its financial support. They also thank Professor H. H. Kuo for making them aware of reference \cite{Kubo}.

\end{document}